\documentclass{article}
\usepackage{amsmath, amssymb, amsfonts, amsthm}
\usepackage{amscd}

\theoremstyle{plain}
\newtheorem{main}{Theorem}

\newtheorem*{TPF}{The Principle of Functoriality}
\newtheorem{theorem}{Theorem}[section]
\newtheorem{lemma}[theorem]{Lemma}
\newtheorem{corollary}[theorem]{Corollary}
\newtheorem{proposition}[theorem]{Proposition}

\begin{document}

\title{Dimension data, and Local versus Global Conjugacy in Reductive Groups}
\author{Song Wang}        
\date{}          
\maketitle

\hspace{20em}To Prof.\ Robert Langlands

\bigskip

\section{Introduction} \label{S:1}

Let $G$ be a reductive algebraic group and $H$ be a connected
semisimple subgroup of $G$. The dimension data consists of the
collection $\{\,{\rm dim} V^{H}\,\}$ with $V^{H}$ denoting the space
of points in $V$ fixed by $H$, and where $(\rho, V)$ runs through
all representations of $G$. One may ask if the dimension data
determine $H$ up to conjugacy or at least isomorphism. In other
words, If $H$ and $H'$ have the same dimension data, are they
isomorphic? If so, are they conjugate in $G$?

\medskip

For the first question, when $G = GL (n, \mathbb{C})$,  M. ~Larsen
and R. ~Pink (\cite{L-P90}) gave an affirmative answer. i.e., if $H$
and $H'$ are two semisimple subgroups of $GL (n, \mathbb{C})$ and
have the same dimension data, then $H$ and $H'$ must be isomorphic.
Moreover, as each algebraic group embeds into $GL (n)$, the answer
is still ``yes'' for the first question when we replace $G$ by an
arbitrary complex reductive algebraic group.

\medskip

For the second question, M. ~Larsen and R. ~Pink also gave an answer
in the same paper for a special case. They proved when $G = GL (n,
\mathbb{C})$, and $H$, $H'$ irreducible (i.e., not contained in any
proper Levi subgroup of $G$), having the same dimension data forces
them to be conjugate in $G$. This fact is used in our article.

\medskip

Our first result is that, the answer is still ``yes'' if $G = SO (2
n + 1)$, $Sp (2 n)$ and $O (N)$ with certain conditions. More
precisely, we have the following:

\medskip

\begin{main}  \label{TM:O}
Let $\rho$ and $\rho'$ be two homomorphisms from a connected
semisimple complex group $H$ into $G = O (N)$, $SO (N)$, or $Sp (N)$
{\rm ( $N = 2 n$ )}, which is naturally embedded into $GL (N)$. And
assume that

\textnormal{(i)} $\rho$ and $\rho'$ are irreducible in
$GL (N)$;

\textnormal{(ii)} $\rho$ and $\rho'$ possess the same dimension
data.

\noindent Then $\rho (H)$ and $\rho' (H)$ differ by an automorphism
of $G$. Moreover, if $G = O (N)$, $Sp (2 n)$ or $SO (2 n + 1)$, this
automorphism is inner so that $\rho$ and $\rho'$ are globally
conjugate in image.
\end{main}

\medskip

However, this is \emph{not} the case in general. In this article
we give a family of pairs $(G, H)$ with $G = SO (2N, \mathbb{C})$,
such that the dimension data does not determine $H$ up to
conjugacy. In addition, our examples are {\it connected, complex
semisimple groups}, so that they give rise to connected analogues
to known examples due to Larsen and Blasius (\cite{Larsen94},
\cite{Larsen96}, \cite{Blasius94}) where $H$ is {\it discrete and
finite}.

\medskip

Our second result, which is the main one of the paper, is the
following:
\begin{main}  \label{TM:A}
Let $H$ be a simple Lie group over $\mathbb{C}$ with its Lie
algebra being one of the following types:
\[
A_{4 n} (n \ge 1), B_{2 n} (n \ge 2), C_{2 n} (n \ge 2),
    E_{6}, E_{8}, F_{4}, G_{2}
\]
and $G = SO (2 N)$ where $2 N = {\rm dim} Lie (H)$. Then there exist
two embeddings $i$ and $i'$\
of $H$ into $G$,
such that their images $i (H)$ and $i'(H)$ are not
conjugate in $G$,
but they possess the same dimension data. In fact,
$i$ and $i'$ are locally conjugate, but
not globally conjugate in image.
\end{main}

Here ``locally conjugate'' means that $i (h)$ and $i' (h)$ are
conjugate in $G$ for each semisimple $h$; ``globally conjugate in
image'' means that $i (H)$ and $i' (H)$ are conjugate in $G$ (cf.\
Section ~\ref{S:2}).

\medskip

The simplest example comes when $H = SO (5)$, $G = SO (10)$, $i =
Ad$, the adjoint representation, or $\Lambda^{2}$, the exterior
square, and $i' = \tau \circ i$ where $\tau$ is some automorphism on
$SO (10)$ which is not inner. In fact, we consider the situation
with the following conditions: $H = Int (\mathfrak{g})$ where
$\mathfrak{g}$ is a simple Lie algebra of even rank, $i$ is the
adjoint representation $Ad$ of $H$ whose image in $GL
(\mathfrak{g})$ is contained in $G =  SO (\mathfrak{g}, \kappa)
\cong SO (2 N)$ where $\kappa$ is the Killing form of
$\mathfrak{g}$, and $i' = \tau \circ i$ where $\tau$ is some
automorphism on $G = SO (\mathfrak{g}, \kappa)$. The list given in
Theorem ~\ref{TM:A} exhausts all possible cases in this situation.

\bigskip

We will prove Theorem ~\ref{TM:O} in Section ~\ref{S:2} and Theorem
~\ref{TM:A} in Section ~\ref{S:3} and ~\ref{S:4}. One may observe
that our construction and proof are still available if $\mathbb{C}$
is replaced by any algebraic closed field of characteristic $0$ such
as $\bar{\mathbb{Q}}$ or $\bar{\mathbb{Q}}_{l}$. In fact, a
potential application is in comparing for two continuous
homomorphisms $\rho_{l}$ and $\rho'_{l}$ from $\mathcal{G}_{k}$ to
$G (\bar{\mathbb{Q}}_{l})$ for $\mathcal{G}_{k}$ the absolute Galois
group for a number field $k$ such that the Zariski closure $H_{l}$,
$H'_{l}$ of ${\rm Im} (\rho_{l})$ and ${\rm Im} (\rho'_{l})$
respectively have same dimension data. By Tchebotarev, the
Frobenious classes ${\rm Frob}_{v}$, for $v$ unramified in
$\rho_{l}$ and $\rho'_{l}$, generate the Galois images
topologically, and local conjugacy implies that for any algebraic
representation $r$: $G \to GL_{N}$, the $L$-factors $L_{v} (s, r
\circ \rho_{l})$ and $L_{v} (s, r \circ \rho'_{l})$ are the same.

\bigskip

We show in Section ~\ref{S:6} that if certain instances of
Langlands' Principle of functoriality are known, then our examples
give rise to cusp forms $\pi$ on certain $SO (2 n) / F$, $F$ a
number field, which appear with multiplicity bigger than one. The
interest in this is that the conjectural Langlands group $H (\pi)$
in such an example will be a connected subgroup of the dual group
$\hat{G}$. Earlier instances of such failures of multiplicity one
like for $SL (n) (n \geq 3)$ (\cite{Blasius94}) because of the local
vs global conjugacy issues, are associated to $H (\pi)$'s which are
disconnected, even finite. For $SL (2)$, one knows multiplicity one
(\cite{Ra2000}), and for $n \geq 3$, E.M. ~Lapid (\cite{Lapid99})
showed that under a Tanakian formalism, the multiplicities are
bounded for each $n$. However, this is not the case for $G_{2}$
(\cite{Gan}) \footnote{For $G_{2}$ in \cite{Gan}, the reason of
failure of multiplicity one is different from that for $SL (n)$.
Indeed, in $G_{2}$, it is expected that local conjugacy implies
global conjugacy. The reason for the failure of multiplicity one in
$G_{2}$ is because the Arthur multiplicity formula gives answers
which are $> 1$. For classical groups, Arthur's formula always gives
$1$ or $0$; yet multiplicity one fails at the times because
local-global conjugacy does not hold.}.

\bigskip

\begin{main}  \label{TM:C}
Let $2 N = {\rm dim} Lie (H)$, where $H$ is a simple Lie group of
one of the following type
\[
A_{4 n} (n \ge 1), B_{2 n} (n \ge 1), C_{2 n} (n \ge 1),
    E_{6}, E_{8}, F_{4}, G_{2}
\]

Assume the following conditions:

\medskip

\textnormal{(1)} The Langlands functoriality from $H$ to $SO (2 N)$
which arises from the adjoint representation $i$ of $H$ holds. Also,
The Langlands functoriality for $i' = \tau \circ i$ also holds,
where $\tau$ is an outer automorphism of $SO (2 N)$.

\medskip

\textnormal{(2)} The Arthur's conjecture, including the Arthur's
multiplicity formula for $SO (2 N)$ holds.

\medskip

\textnormal{(3)} There is a cuspidal automorphic representation
$\pi_{0}$ of $H$, which is associated to a global $l$-adic Galois
representation $\rho_{l}$ with ${\rm Im} (\rho_{l})$ Zariski dense
in the dual group $\hat{H} (\bar{\mathbb{Q}}_{l})$.

\medskip

Then the multiplicity one fails for $G = SO (2 N)$. Moreover, there
are two nearly equivalent cuspidal representations $\pi$ and $\pi'$,
not of finite Galois type, which cannot be distinguished by their
incomplete $L$-functions, i.e., for a finite set $S$ of places,
\[
L^{S} (s, \pi, r) = L^{S} (s, \pi', r)
\]
for any algebraic representation $r$ of $\hat{G}$; yet they don't
occur in the same constituent and hence give multiplicity $> 1$.

\end{main}

\medskip

For a proof, see Section ~\ref{S:6}.

\bigskip

This result was suggested by Langland's paper ``Beyond Endoscopy''
(\cite{Langlands2002}), and this Theorem gives nontrivial evidence
for his prediction.

\bigskip

In a sequel we plan to investigate the instances of functoriality
when $F$ is a function field over a finite field, by appealing to
Lafforgue's work as well as the recent forward and backward lifting
results of generic forms from the classical groups to $GL (n)$. The
transfer in the forward direction was done for number fields by
Cogdell, Kim, Piatetski-Shapiro and Shahidi in \cite{CoKPSS} (and
the papers to follow), and very recently, progress on the analogue
for function fields over finite fields was done by L.A.~ Lomeli in
his Purdue thesis. The failure of multiplicity one for suitable $SO
(2 n)$ comes from the outer automorphism $\tau$ of $SO (2 n)$ which
arises from a conjugation in $O (2 n)$. The fact that cusp forms on
$SO (2 n)$ which are moved by the outer automorphism will appear
with multiplicity is one of the key observations of this work and it
points to a need to fine tune Arthur's original multiplicity
conjectures. It raises the following question: Can one use the
twisted trace formula relative to $\tau$ and isolate the cusp forms
with multiplicity one. We hope to investigate this in another
sequel.

\bigskip


I would like to express our thanks to R.P.~ Langlands for suggesting
this problem, and his helpful advice at the IAS, and his positive
remarks on our results. Also I am grateful to D. ~Ramakrishnan for
his continued interest in this problem and for helpful discussions
over the past years. I also thank Wee Tek Gan for reading the first
draft of this paper carefully and for making suggestions for
improvement. I also benefit from his helpful comments on Arthur's
conjecture and multiplicity formula during revision of that earlier
draft. Moreover, I want to acknowledge the inspiration I got from
reading the paper \cite{L-P90} of Larsen and Pink, and the works of
Arthur on the multiplicities of cusp forms on reductive groups.
Finally, the author acknowledges the support from the NSF Grant
DMS--9729992 while being a member of the IAS, Princeton, during
2001--2002.

\bigskip

\section{Local Conjugacy and Global Conjugacy}
\label{S:2}

Let $H$ and $H'$ be two semisimple subgroups of $G$. We say that
they are \emph{\textbf{locally conjugate}} if there is an
isomorphism $\tau$ from $H$ to $H'$ such that $\tau (h)$ and $h$ are
conjugate in $G$ for each semisimple $h \in H$; $H$ and $H'$ are
\emph{\textbf{globally conjugate}} if they differ by a conjugation,
i.e., an inner automorphism of $G$. Moreover, let $i$ and $i'$ be
two embeddings of $H$ into $G$. We say that $i$ and $i'$ are
\emph{\textbf{locally conjugate}} if $i (h)$ and $i' (h)$ are
conjugate in $G$ for each semisimple $h \in H$; $i$ and $i'$ are
\emph{\textbf{globally conjugate in image}} if $i (H)$ and $i' (H)$
are globally conjugate in $G$, i.e., $i' (H) = c (i (H))$ for some
inner automorphism $c$ of $G$; They are \emph{\textbf{globally
conjugate}} if $i' = c \circ i$ for some inner automorphism $c$ of
$G$.

\bigskip

\emph{Remark}: If $i, i'$: $H \hookrightarrow G$ are globally
conjugate in image, then there exists an automorphism $\tau$ of $H$,
and an element $t \in G$ such that, for each $h \in H$, $i' (\tau
(h)) = t^{-1} i (h) t$.

\bigskip

\begin{proposition} \label{T:201}
If $H$ and $H'$ are locally conjugate in $G$, then they have the
same dimension data.
\end{proposition}

\medskip

\emph{Proof. }  Let $i$ be the identity map on $H$ and $i'$ be the
isomorphism from $H$ to $H'$ such that $i' (h)$ and $i (h)$ are
conjugate for each $h \in K (H)$ where $K (H)$ is a maximal
compact subgroup of $H$. This is possible since all elements on a
maximal compact subgroup are semisimple. Then for each
representation $(\sigma, V)$ of $G$, we have

\begin{align}
{\rm dim V^{H}} &= {\rm mult} (1, \sigma \circ i) \notag \\
&= \frac{1}{{\rm Vol} (K (H))} \int_{K (H)} {\rm Tr} \sigma
(h)\,d\mu \notag
\end{align}

\noindent and similarly,

\begin{align}
{\rm dim V^{H'}} &= {\rm mult} (1, \sigma \circ i') \notag \\
&= \frac{1}{{\rm Vol} (K (H))} \int_{K (H)} {\rm Tr} (\sigma \circ
i') (h) \,d\mu \notag
\end{align}

\medskip

Thus to prove that $H$ and $H'$ have the same dimension data, it
suffices to prove the following claim: ${\rm Tr} \sigma (h) = {\rm
Tr} (\sigma \circ i') (h)$ for any $h$ in $K(H)$. As $H$ and $H'$
are locally conjugate, $i (h) = h$ and $i' (h)$ are conjugate in
$G$ for each $h$ in $K (H)$. So $i' (h) = \beta^{-1} h \beta$ for
some $\beta$ in $G$. So, $(\sigma \circ i') (h) = \sigma
(\beta)^{-1}\, \sigma (h)\, \sigma (\beta)$. Hence the claim and
the theorem.

\qed

Remark: This theorem still applies for $G$ over
$\bar{\mathbb{Q}}_{l}$,
the algebraic closure of $\mathbb{Q}_{l}$. The reason is that,
for any algebraic subgroup $H$, and any algebraic
representation $(\sigma, V_{K'})$ of $G$,
\[
{\rm dim}_{\bar{\mathbb{Q}}_{l}} V_{\bar{\mathbb{Q}}_{l}}^{H}
= {\rm dim}_{K'} V_{K'}^{H}
= {\rm dim}_{\mathbb{C}} V_{\mathbb{C}}^{H}
\]
where $K'$ is the number field where $i$ and $\sigma$ split.

\bigskip

Next, we prove Theorem ~\ref{TM:O}.

\medskip

Let $\rho$ be an irreducible selfdual representation of a complex
semisimple group $H$ into $GL (V)$. Then $\rho$ is either
orthogonal or symplectic, and in either case, the image of $\rho$
must fix some symmetric or alternating nondegenerate bilinear form
$\omega$ on $V$. In fact, such form must be unique up to scalar.

\medskip

\begin{lemma} \label{T:202}
Let $\rho$ be an irreducible selfdual finite dimensional
representation of a connected semisimple group $H$. Then There is a
unique nondegenerate bilinear form $\omega$ up to a scalar that $H$
preserves. i.e., $\rho (H) (\omega) = \omega$.
\end{lemma}

\medskip

\emph{Proof.  } Note that each nondegenerate bilinear form, either
symmetric or alternating, up to a scalar, that $H$ fixes (or fixes
up to a scalar factor) corresponds a trivial (or a one dimensional)
constituent of $\rho \otimes \rho$. Hence, it suffices to prove
that, $\rho \otimes \rho$ possesses exactly one one-dimensional
constituent. This is the case since following: First, $H$ is
semisimple and selfdual. and thus $\rho$ is equivalent to
$\tilde{\rho}$ which is the contragradient of $\rho$. Moreover, the
multiplicity of one occurred in $\rho \otimes \tilde{\rho}$ is
exactly the square sum of the multiplicities of the irreducible
constituents of $\rho \otimes \tilde{\rho}$, which is $1$ if $\rho$
is irreducible.

\qed

\bigskip

\begin{lemma} \label{T:203}
Let $\rho$ and $\rho'$ be two irreducible representations of a
connected complex semisimple group $H$ into $GL (V) = GL (N,
\mathbb{C})$. Assume that

\textnormal{(i)} $ \rho (H)$ and $\rho' (H)$ fixes the same
nondegenerate bilinear form $\omega$, either symmetric or
alternating;

\textnormal{(ii)} $\rho$ and $\rho'$ are globally conjugate
in $GL (V)$.

Then they are globally conjugate in $O (V, \omega)$,
i.e., there is a $t \in GL (V)$ such that,
$\rho' (g) = t^{-1} \rho (g) t$ for each $g$ in $H$,
and $t$ fixes $\omega$.
\end{lemma}

\medskip

\emph{Proof. } Let $t$ be an element in $GL (V)$ such that $\rho'
(g) = t^{-1} \rho (g) t$ for each $g$ in $H$. Then we claim that
some scalar multiple of $t$, namely, $t_{1} = c t$, fixes $\omega$.
Granting this, we have $t_{1} \in O (V, \omega)$, and $\rho' (g) =
t_{1}^{-1} \rho (g) t_{1}$ for each $g$ in $H$. Hence the lemma.

\medskip

Now we prove the claim. Check that,
$H$ fixes $t \omega$. In fact, for each $h \in H$,
\[
\rho (h) (t \omega) = (\rho (h) t) (\omega) =
(t \rho' (h)) (\omega) = t (\rho' (h) (\omega)) = t \omega
\]
So by Lemma ~\ref{T:202}, together with the assumption
that $\rho$ is irreducible, $t \omega = c_{1} \omega$
for some nonzero $c_{1}$. Thus the claim.

\qed

\bigskip

\textbf{\emph{Proof of Theorem ~\ref{TM:O}}}.

\medskip

First we treat the case when $G = O (N)$ or $Sp (N)$. Then $G = O
(N, \omega)$ for some nondegenerate bilinear form $\omega$, either
symmetric or symplectic on $V$ where $V$ is the $N$--dimensional
space where $G \subset GL (N)$ acts.

\medskip

Note that the dimension data of $(H, G)$ are a sub-collection of the
dimension data of $(H, GL (N))$. Then (i), (ii), plus the
Larsen--Pink imply that $\rho$ and $\rho'$ are globally conjugate in
image in $GL (N)$. Hence there is some $\varphi \in {\rm Aut} (H)$
such that $\rho \circ \varphi$ and $\rho'$ are globally conjugate in
$GL (N)$. Note that $\rho \circ \varphi$ is also irreducible. Hence
by Lemma ~\ref{T:203}, $\rho \circ \varphi$ and $\rho'$ are globally
conjugate in $O (N)$. Hence the theorem.

\medskip

For $G = SO (2 n + 1)$, we also have the theorem as $O (2 n + 1) =
\{\pm I\} SO (2 n + 1)$.

\medskip

\qed

\bigskip

\emph{Remark:} One expects that, if $G$ does not allow outer
automorphisms, then the second question raised in the introduction
of this article is affirmative, with only the assumption that, $H$
or $H'$ are not contained in any proper parabolic or Levi subgroups
of $G$.

\bigskip

\section{Proof of Theorem B, Part I} \label{S:3}

Let $\mathfrak{g}$ be a simple Lie algebra with the Killing form
$\kappa$ defined as $\kappa (X, Y) = {\rm Tr} (ad (X) ad (Y))$. It
is clear that $\kappa$ is nondegenerate by the semisimplicity. Let
$H$ be the adjoint semisimple Lie group with Lie algebra
$\mathfrak{g}$. Here $adjoint$ means $Z (H)$ is trivial. Such $H$
exists and is unique by the Lie theory. In fact, $H \simeq {\rm
Int} (\mathfrak{g})$. Hence, the adjoint representation of $H$
gives rise to an irreducible representation of $H$ in the
underlying space of $\mathfrak{g}$:

\[
i: H \hookrightarrow Gl (\mathfrak{g})
\]

\bigskip

\begin{lemma}  \label{T:301}
$H$ preserves the Killing form. So that, the image of $H$ lies in
$SO (\mathfrak{g}, \kappa)$.
\end{lemma}

\medskip

Here for each linear space $V$ and a nondegenerate symmetric
bilinear form $\omega$, $SO (V, \omega)$ consists of linear
transforms that preserves $\omega$. It is clear that $SO (V,
\omega)$ is isomorphic to $SO ({\rm dim} V)$.

\bigskip

\emph{Proof}. For each $X$, $Y$ and $Z$ in $\mathfrak{g}$, $c \in
{\rm Int} (\mathfrak{g})$,

\begin{align}
{\rm ad} (c(X))& {\rm ad} (c(Y))\, (Z) = [c(X),\, [c(Y),\, Z]]
\notag \\
&= [c(X),\, c([Y,\, c^{-1}(Z)])] = c([X,\, [Y,\, c^{-1} (Z)]])
\notag \\
&= c({\rm ad} (X) {\rm ad} (Y)\, c^{-1} (Z)) \notag
\end{align}

Hence ${\rm ad} (c(X)) {\rm ad} (c(Y)) = c \circ {\rm ad} (X) {\rm
ad} (Y) \circ c^{-1}$ is conjugate to ${\rm ad} (X) {\rm ad} (Y)$.
So they must have the same trace.

\medskip

Furthermore, $H$ is connected and so is its image. Therefore, the
image must lie in $SO (\mathfrak{g}, \kappa)$, which is the
identity connected component of $O (\mathfrak{g}, \kappa)$. Done.

\qed

\bigskip

Now we make several assumptions to help to explain our examples.

\medskip

\underline{Assumption (A)}: the rank of $\mathfrak{g}$ is even.

\medskip

So by Lemma ~\ref{T:301} under Assumption (A), $i$ is orthogonal,
and its image lies in $SO (\mathfrak{g}, \kappa)$ which is
isomorphic to $SO (2 N, \mathbb{C})$ where $2 N$ is the dimension of
$\mathfrak{g}$.

\bigskip

Let $\tau$ be any ``odd'' element in $O (\mathfrak{g}, \kappa)$,
i.e., $\tau$ preserves $\kappa$ but does not lie in the identity
component $SO (\mathfrak{g}, \kappa)$. Also, let $i' = {\rm Conj}
(\tau) \circ i$, i.e., $i' (x) = \tau i (x) \tau^{-1}$.

\medskip

Let $\Phi$ be the root system of $H$ and $\Delta$ be a base. We
first assume the following:

\medskip

\underline{Assumption (B)}: ${\rm Aut} (\Phi) = {\rm Inn} (\Phi)$.

\medskip

Assumption (B) is equivalent to that, the only isometry of the
Dynkin diagram is the identity map. Under such assumption, the Lie
automorphism of $H$ must be a conjugacy by some element in $G$. The
simple Lie algebras satisfying Assumption (B) are $A_{1}$, $B_{n} \,
(n \geq 2)$, $C_{n} \, (n \geq 2)$, $E_{7}$, $E_{8}$, $F_{4}$ and
$G_{2}$.

\medskip

In next section, we will deal with the case when $\Phi$
allows a outer automorphism.

\medskip

\begin{proposition} \label{T:302}
With the Assumption \textnormal{(A)} and \textnormal{(B)}, $i$ and
$i'$ are locally conjugate, but not globally conjugate in image. The
only cases that apply are $\mathfrak{g} = B_{2 n}$, $C_{2 n}$,
$E_{8}$, $F_{4}$ and $G_{2}$.
\end{proposition}

\bigskip

\emph{Remark:} So the simplest example should be $\mathfrak{g} =
B_{2} = C_{2}$. In this example, $H$ is $SO (5)$, and the adjoint
representation of $H$ gives rise to the exterior square from $SO
(5)$ to $SO (10)$.

\bigskip

We need two lemmas for the proof.

\medskip

\begin{lemma} \label{T:303}
Let $A$ and $B$ are two diagonalizable matrices in $SO (2 N)$, and
suppose that they are conjugate in $O (2 N)$, and $1$ occurs as an
eigenvalue for either $A$ or $B$. Then $A$ and $B$ are conjugate in
$SO (2 N)$.
\end{lemma}

\emph{Proof}. First, we claim that $A$ commutes with an ``odd''
element in $O (2 N)$, i.e., an element in $O (2 N)$ but not $SO (2
N)$.

\medskip

We choose an appropriate coordinate system so that the invariant
quadratic form on $O (2 N, \mathbb{C})$ is $\omega = \sum x_{i}
y_{i}$. So, with respect to a basis \newline $\{ e_{1}, f_{1},
e_{2}, f_{2}, \ldots, e_{n}, f_{n} \}$, the matrix representing
$\omega$ is ${\rm Diag} (P, P, \ldots, P)$ where $P =
\begin{pmatrix} 0 & 1
\\ 1 & 0 \end{pmatrix}$. By conjugation in $SO (2 N)$ and without loss of
generality, we may also assume $A$ is diagonal, i.e., of the form
${\rm diag} (a_{1}, a_{1}^{-1}, a_{2}, a_{2}^{-1}, \ldots, a_{n},
a_{n}^{-1})$. So if say $a_{1} = 1$, then the centralizer of $A$
contains $O (2) \times Id$ where $O (2)$ is the orthogonal group on
$\langle e_{1}, f_{1} \rangle$ and $Id$ is the identity map on the
orthogonal complement $\langle e_{2}, f_{2}, \ldots, e_{n}, f_{n}
\rangle$. Hence the claim.

\medskip

Now assume that $B = C^{-1} A C$. If $C \in SO (2 N)$ we are done.
Otherwise, $C$ is an ``odd'' element. Hence,

\[ B = C^{-1} A C = C^{-1} Q^{-1} A Q C = {(Q C)}^{-1} A (Q C)
\]
for some ``odd'' element $Q$ which commutes with $A$. Thus $Q C \in
SO (2 N)$, and hence $A$ and $B$ must be also conjugate in $SO (2
n)$.

\qed

\bigskip

\begin{lemma} \label{T:304}
Let $\rho:$ $H \to GL (n, \mathbb{C})$ be an irreducible orthogonal
representation. Then the centralizer of $H$ in $O (n, \mathbb{C})$
is $\{\, \pm I\,\}$.
\end{lemma}

\emph{Proof}. As $\rho$ is irreducible, $C_{GL (n, \mathbb{C})}
(\rho (H)) = C^{*}$ from the Schur's Lemma. Then $C_{O (n,
\mathbb{C})} (\rho (H)) = C_{GL (n, \mathbb{C})} \cap O (n,
\mathbb{C}) = C^{*} \cap O (n, \mathbb{C}) = \{ \pm I \}$.

\qed

\bigskip

\textbf{\emph{Proof of Proposition ~\ref{T:302}}}.

\medskip

{\emph{Step 1.}  $i$ and $i'$ are
locally conjugate.}

\medskip

For each $h_{0}$ semisimple, $h = i (h_{0})$ and $h' = i' (h_{0})
= \tau h \tau^{-1}$, where $\tau \in O (\mathfrak{g}, \kappa)
\simeq O (2 N)$. Hence $h$ and $h'$ are conjugate in $O
(\mathfrak{g}, \kappa)$. Moreover, they are semisimple also as $i$
and $i'$, and hence each of them must be contained in some complex
torus, so that $h$ fixes some Cartan subalgebra in $\mathfrak{g}$.
This shows that $h$ has eigenvalue $1$ in $SO (\mathfrak{g},
\kappa)$. Furthermore the multiplicity of $1$ is at least the
dimension of this Cartan Subalgebra, namely, the rank of
$\mathfrak{g}$ which is at least $2$ from Assumption (A) which
says that the rank of $\mathfrak{g}$ is even. Thus Lemma
~\ref{T:303} applies, and $h$ and $h'$ must be conjugate in $SO
(\mathfrak{g}, \kappa)$. By the arbitrary choice of $h_{0}$, $i$
and $i'$ are locally conjugate.

\medskip

{\emph{Step 2.}  $i$ and $i'$ are
not globally conjugate in image.}

\medskip

Now we need to use Assumption (B). Assume the contrary, that $i$
and $i'$ are globally conjugate in image, and moreover, that the
conjugation by $\beta \in SO (\mathfrak{g}, \kappa)$ sends $i'
(H)$ to $i (H)$, i.e., $i (H) = \beta i' (H) \beta^{-1}$. Thus the
conjugation by $\beta \tau$ restricted to $i (H)$ is an
automorphism on $i (H) \simeq H$. From Assumption (B), all Lie
automorphisms on $i (H) \simeq H$ must be inner. Then this
automorphism is in fact the conjugation by some element $\gamma$
in $i (H) \subset SO (\mathfrak{g}, \kappa)$. Hence $\gamma^{-1}
\beta \tau$ must lie in the centralizer of $i (H)$ in $O
(\mathfrak{g}, \kappa) \simeq O (2 N)$ since the conjugation by
$\gamma^{-1} \beta \tau$ fixes $i (H)$ at all. As $\mathfrak{g}$
is simple, then the adjoint representation of $H$ in
$\mathfrak{g}$ must be irreducible. Hence by Lemma ~\ref{T:304},
$\gamma^{-1} \beta \tau$ must be $\pm I$ which is obviously an
element in $SO (\mathfrak{g}, \kappa)$. However, $\tau$ is
``odd'', and $\beta$ and $\gamma$ are ``even'' as they lie in $SO
(\mathfrak{g}, \kappa)$. Then $\gamma^{-1} \beta \tau$ must be
``odd''. This gives a contradiction.

\medskip

Hence $i$ and $i'$ are not globally
conjugate in image.

\bigskip

Then the first part of this theorem follows. The second part is
clear as the only connected Dynkin diagrams of even rank without any
outer automorphisms are $B_{2 n}$, $C_{2 n}$, $E_{8}$, $F_{4}$ and
$G_{2}$.

\qed

\bigskip

\emph{Remark: } When Proposition ~\ref{T:302} applies, $i (H)$ and
$i' (H)$ have the same dimension data. But they are not
necessarily conjugate in $G = SO (\mathfrak{g}, \kappa)$.

\bigskip

\section{Proof of Theorem B, Part II} \label{S:4}

Consider the adjoint representations of simple adjoint groups. What
happens if the root system $\Phi$ admits a nontrivial outer
automorphism? For example, consider the type $A_{n}$, $D_{n} (n \geq
3)$, $E_{6}$. We will prove, if $\mathfrak{g}$ is $A_{4 n}$ or
$E_{6}$, we still have the same conclusion, i.e, $i$ and $i'$ are
locally conjugate but not globally conjugate in image.

\medskip

In fact, it is enough for us to replace Assumption (B) by the
following weaker assumption.

\medskip

\underline{Assumption (C)}: All elements of ${\rm Aut}(H)$ are
``even'' in the adjoint representation.

\bigskip

First, we claim that, given an adjoint representation $i: H = {\rm
Int}(\mathfrak{g}) \to GL(\mathfrak{g})$ where $\mathfrak{g}$ is
an simple Lie algebra of even degree, if $\tau_{1}$ is an
automorphism of $H$, which induces an automorphism $\tau'_{1}$ of
$\mathfrak{g}$, then $\tau_{1}$ is the restriction of ${\rm Conj}
(\tau'_{1})$ on $GL (\mathfrak{g})$ to $H$, i.e.,
\[
    {\rm Ad} (\tau_{1} (h)) (X) = \tau'_{1} \circ {\rm Ad} (h)
     \circ {\tau'}_{1}^{-1} (X)
\]
for each $h \in H$ and $X \in \mathfrak{g}$.

\bigskip

\noindent {\bf Definition} Assume that the rank of $\mathfrak{g}$
is even. $\tau_{1}$ is said to be ``even'' if $\tau'_{1}$, viewed
as a linear transformation, has determinant $1$, and ``odd'' if
$\tau'_{1}$ has determinant $-1$.

\bigskip

\emph{Remark. } If $\tau_{1}$ is inner, then it must be ``even''
since $\tau'_{1}$ lies in the image of $H$ in $GL (\mathfrak{g})$,
hence in $SO (\mathfrak{g}, \kappa)$ as $H$ is connected, thus
having determinant $1$.

\bigskip

\begin{lemma} \label{T:401}
For each $\tau_{1} \in {\rm Aut} (H)$, $\tau'_{1}$ lies in $O
(\mathfrak{g}, \kappa)$. Hence $\tau_{1}$ is ``even'' if and only if
$\tau'_{1}$ lies in $SO (\mathfrak{g}, \kappa)$. In particular, all
the inner automorphisms of $H$ are ``even''.
\end{lemma}

\medskip

\emph{Proof}.  The idea is almost the same as in the proof of
Lemma ~\ref{T:301}. For each $X$, $Y$ and $Z$ in $\mathfrak{g}$,
we have
\begin{align}
{\rm ad} (\tau'_{1} (X))& {\rm ad} (\tau'_{1}(Y))
\, (Z) = [\tau'_{1} (X),\, [\tau'_{1} (Y),\, Z]]
\notag \\
&= [\tau'_{1} (X),\, \tau'_{1} ([Y,\, {\tau'_{1}}^{-1}(Z)])]
 = \tau'_{1} ([X,\, [Y,\, {\tau'_{1}}^{-1} (Z)]])
\notag \\
&= \tau'_{1} ({\rm ad} (X) {\rm ad} (Y)\, {\tau'_{1}}^{-1} (Z)) \notag
\end{align}

Here we used the fact that $\tau'_{1}$ preserves the Lie brackets
as it come from ${\rm Aut} (H)$ and hence a Lie algebra
automorphism.

\medskip

Hence ${\rm ad} (\tau'_{1} (X)) {\rm ad} (\tau'_{1}(Y))$
and ${\rm ad} (X)  {\rm ad} (Y)$ have the same trace. Hence
$\tau'_{1}$ preserves the Killing form of $\mathfrak{g}$.
The rest assertions are now clear.

\medskip

\qed

\bigskip

Recall that $i$ denotes the adjoint representation
of $H$, and $i' = {\rm Conj} (\tau) \circ i$.
Here $H$ is an adjoint group with simple Lie algebra $\mathfrak{g}$.

\medskip

\begin{proposition}  \label{T:402}
With the Assumption \textnormal{(A)} and \textnormal{(C)}, $i$ and
$i'$ are locally conjugate, but not globally conjugate in image.
\end{proposition}

\bigskip

\emph{\textbf{Proof of Proposition ~\ref{T:402}}}.

\medskip

Recall that in the proof of Proposition ~\ref{T:302}, Assumption
(B) is used only in Step 2. Hence the local conjugacy of $i$ and
$i'$ is clear when we replace Assumption (B) with (C). Thus it
suffices to redo Step 2.

\medskip

{\emph{Step 2$'$.}  $i$ and $i'$ are not globally conjugate in
image.}

\medskip

Now we need to use Assumption (C) instead of Assumption (B).
Recall that $i' (H) = \tau^{-1} i (H) \tau$. Now again assume the
contrary, i.e., $i$ and $i'$ are globally conjugate in image, and
again that $i (H) = \beta i' (H) \beta^{-1}$ for $\beta \in SO
(\mathfrak{g}, \kappa)$. Then the conjugation by $\beta \tau$
restricted to $i (H) \sim H$ is an automorphism on $H$. From
Assumption (C), and Lemma ~\ref{T:401}, it induces the conjugation
by $\gamma \in SO (\mathfrak{g}, \kappa)$. Hence $\gamma^{-1}
\beta \tau$ must lie in the centralizer of $i (H)$ in $O
(\mathfrak{g}, \kappa) \simeq O (2 N)$ since the conjugation by
$\gamma^{-1} \beta \tau$ fixes $i (H)$ at all. As $\mathfrak{g}$
is simple, then the adjoint representation of $H$ in
$\mathfrak{g}$ must be irreducible. Hence by Lemma ~\ref{T:304},
$\gamma^{-1} \beta \tau$ must be $\pm I$ which is obviously an
element in $SO (\mathfrak{g}, \kappa)$. However, $\tau$ is
``odd'', and $\beta$ and $\gamma$ are ``even'' as they lie in $SO
(\mathfrak{g}, \kappa)$. Then $\gamma^{-1} \beta \tau$ must be
``odd''. This gives a contradiction also.

\medskip

Hence $i$ and $i'$ are not globally
conjugate in image.

\qed

\bigskip

Now we want to refine further and prove that Assumption (C) can in
fact be replaced by the following equivalent assumption.

\medskip

\underline{Assumption (C$^{\prime}$)}: All automorphisms of $\Phi$
that preserves a base is an even permutation on the base. i.e.,\
All automorphisms of the Dynkin diagram are even permutations of
vertices.

\bigskip

In fact, since each automorphism of $H$ is a product of an inner
automorphism which is ``even'' by Lemma ~\ref{T:401}, and an
automorphism that preserves a maximal torus $T$ and a base
$\Delta$ of the root system, we may focus on those $\tau_{1}$ that
preserves $T$ and $\Delta$. Denote $\mathfrak{t}$ be the Cartan
subalgebra corresponding to $T$. Then $\tau'_{1}$, induced from
$\tau_{1}$ as a Lie algebra automorphism on $\mathfrak{g}$,
preserves $\mathfrak{t}$. For each positive root $\beta$, let
\[
V (\beta) = \mathfrak{g}_{\beta} + \mathfrak{g}_{-
\beta}
\]
where $\mathfrak{g}_{\pm \beta}$ is the root space of $\pm \beta$
respectively. Then $\mathfrak{g}$ is a direct sum of $V = \oplus V
(\beta)$ and $\mathfrak{t}$.

\medskip

Hence, as $\tau_{1}$ permutes positive roots, $\tau_{1}$ permutes
$V (\beta)$ for all positive roots $\beta$. In fact, $\tau'_{1}$
sends $\mathfrak{g}_{\beta}$ to $\mathfrak{g}_{\tau_{1} (\beta)}$,
and $\mathfrak{g}_{-\beta}$ to $\mathfrak{g}_{\tau_{1} (-\beta)}$.
Thus $\tau'_{1}$ restricted to $V$ is ``even'', i.e, it has
determinant $1$, and the sign of $\tau'_{1}$ is hence the same as
the sign of $\tau'_{1}$ restricted to $\mathfrak{t}$.

\medskip

Moreover, as $\tau'_{1}$ sends $H_{\beta}$ to $H_{\tau_{1} (\beta)}$
for each $\beta \in \Delta$, where $H_{\beta}$ is the co-root
associated to $\beta$ in $\mathfrak{t}$, the sign of the restriction
of $\tau'_{1}$ to $\mathfrak{t}$ is the same as the sign of
permutation on $\Delta$ by $\tau_{1}$, and it is the same as the
sign of the automorphism of Dynkin graph induced by $\tau_{1}$.
Thus, we have proved the following:

\medskip

\begin{lemma}  \label{T:403}
Assumption \textnormal{(C$^{\prime}$)} is equivalent to Assumption
\textnormal{(C)}.
\end{lemma}

\bigskip

\begin{corollary} \label{T:404}
The only additional cases where Proposition ~\ref{T:402} applies
are $\mathfrak{g} = A_{4 n}$ and $E_{6}$.
\end{corollary}

\bigskip

\emph{Proof}. Let $\mathfrak{g}$ be a simple Lie algebra of even
rank, and assuming that $\mathfrak{g}$ admits an outer automorphism.
Then $\mathfrak{g}$ is of type $A_{2 n}$, $D_{2 n}$ or $E_{6}$.
Suppose in addition that Assumption (C) holds, i.e., all
automorphisms of the Dynkin graph are even permutations. For $A_{2
n}$, the Dynkin graph allows only one outer automorphism, which is a
product of $n$ transpositions on the vertices. Thus Assumption
(C$^{\prime}$) will rule out $A_{2 n}$ when $n$ is odd, and keep
$A_{4 n}$. $D_{2 n}$ is ruled out also as the Dynkin diagram
obviously admits an automorphism which swaps two nodes and keeps
others. $E_{6}$ remains as the only outer automorphism of the Dynkin
graph is a product of two transpositions on the vertices.

\medskip

Hence, by Lemma ~\ref{T:403}, these are the only cases Proposition
~\ref{T:402} applies.

\qed  \bigskip

So Theorem ~\ref{TM:A} follows by combining of Proposition
~\ref{T:302} and Proposition ~\ref{T:402}. When Assumption (C), or
equivalently Assumption (C$^{\prime}$) fails, $\Phi$ and hence
$\mathfrak{g}$ will allows an ``odd'' automorphism, so that any
two subgroups of $SO (2 N)$ conjugate in $O (2 N)$ must be also
conjugate in $SO (2 N)$. Hence $i$ and $i'$ are also globally
conjugate. This explains why the list given in Theorem ~\ref{TM:A}
exhausts all possibilities in this situation.

\bigskip

\section{Preliminaries on Langlands Functoriality \\ and Multiplicity One}
\label{S:5}

In this part, let us recall some preliminary facts on Langlands
Functoriality and the multiplicity one. The experts can skip this
section and go directly to Section ~\ref{S:6}.

\medskip

Fixing a ground field $F$, a local or global field, let $G$ be an
algebraic group, $\mathcal{A} (G)$ the set of automorphic
representations of $G (\mathbb{A}_{F})$ where $\mathbb{A}_{F}$
denotes the adele ring of $F$, and $\mathcal{A}_{0} (F)$ denote the
set of cuspidal automorphic representations of $G (\mathbb{A}_{F})$.
Also, ${}^{L}G = \hat{G} \rtimes \mathcal{W}_{F}$ denotes the
Langlands dual group of $G$, and ${}^{L}G^{\circ} = \hat{G}$ is the
connected component of ${}^{L}G$. For details, see \cite{Bo} and
\cite{Cogdell2003-1}. For the known instances of the local
Langlands, see also \cite{Langlands89} (archimedean case),
\cite{LRS93} (non-archimedean case, positive characteristic),
\cite{HaT2001} (Non-archimedean case, characteristic $0$), and
\cite{Jiang-So2003}, \cite{Jiang-So2002} (for $SO_{2 n + 1}$).

\bigskip

Now, we state the global Langlands conjecture, which is not proved
in general. Now let $F$ be a global field. We need conjectured
$\mathcal{L}_{F}$, the Langlands group of $F$ taking the role of
local Weil or Weil--Delign groups (for the description, see
\cite{LL79}, \cite{Arthur05}). When $\pi$ is of Galois type, then
the global parametrization $\phi$ can be taken as a Galois
representation. Moreover, when we use $l$--adic representations
instead of $\mathbb{C}$--representations, the global Weil group
$\mathcal{W}_{F}$ is believed to serve for this role.

\bigskip

If the local Langlands for each local group $G_{v} = G (F_{v})$ is
assumed, then attaching to $\pi = \bigotimes' \pi_{v}$ we give a
collection of local admissible representations $\{\phi_{v}\}$ of
local parameters $\phi_{v} = \phi_{\pi_{v}}$: $\mathcal{W}'_{F_{v}}
\to {}^{L} G (F_{v})$. Let $\iota_{v}$ be the natural embedding of
${}^{L}G_{v} = \hat{G} \rtimes \mathcal{W}_{F_{v}}$ into ${}^{L}G$.
Then, one gets a collection $\{\iota_{v} \circ \phi_{v}\}$ of local
parameters $\iota_{v} \circ \phi_{v}$: $\mathcal{W}'_{F_{v}} \to
{}^{L}G$.

\medskip

For any finite set of places $S$ of $F$, we can define a family of
global (incomplete) $L$--functions of $\pi$ as
\[
L^{S} (s, \pi, r) = \prod_{v \notin S} L (s, \pi_{v}, r_{v}) =
\prod_{v \notin S} L (s, r \circ \iota_{v} \circ \phi_{v})
\]
when $r$ runs through all representations from ${}^{L}G \to GL_{N}
(\mathbb{C})$ and $r_{v} = r \circ \iota_{v}$.

\bigskip

\begin{TPF} \textnormal{(cf.\ \cite{Langlands2002})}

If $k$ is a local or global field,
$H$ and $G$ be connected reductive $k$--groups
with $G$ quasisplit, then to each $L$--homomorphism
$u$: ${}^{L}H \to {}^{L}G$ there is associated
a transfer/lifting of admissible or automorphic
representations of $H$ to admissible or automorphic
representations of $G$.
\end{TPF}

\medskip

Let $\pi = \bigotimes \pi_{v}$ and $\pi' = \bigotimes \pi'_{v}$ be
two automorphic representations of $G$, and $\{\rho_{v}\}$ and
$\{\rho'_{v}\}$ are the families of local Weil representations to
${}^{L}G$ associate to $\pi$ and $\pi'$ respectively. We say that,
$\pi$ and $\pi'$ are ``locally conjugate'', if for almost all $v$,
$\pi_{v} \cong \pi'_{v}$, i.e., $\pi$ and $\pi'$ are nearly
equivalent, or equivalently, $\rho_{v}$ and $\rho'_{v}$ are
``locally conjugate'' in ${}^{L}G$ for almost all $v$. Since for
almost all $v$ the local $L$-factor $\rho_{v}$ is unramified, and
uniquely determined by the image of the Frobenious elements, then
$\pi$ and $\pi'$ cannot be distinguished by their incomplete
$L$--functions if they are locally conjugate.

\bigskip

Let $\pi = \bigotimes \pi_{v}$ be an automorphic representation of
$G (\mathbb{A}_{F})$. $m (\pi)$, the multiplicity of $\pi$ is
defined as the multiplicity of $\pi$ occurred in \linebreak $L_{\rm
cusp}^{2} (G (F) \backslash G (\mathbb{A}_{F}))$, and it is positive
if and only if $\pi$ is cuspidal. Fixing a conjectural global
parameter $\rho$ of $\pi$, $m (\mathcal{L} (\pi)) = m (\mathcal{L}
(\pi, \rho))$, The multiplicity of an (stable) $L$-packet of $\pi$
(which is singleton if $G = GL (n)$), is defined as the multiplicity
of any $\pi'$ in $\mathcal{L} (\pi)$ which is the set of cuspidal
automorphic representations $\pi'$ associated to the same global
parameter $\rho$ as $\pi$, and it is the number which Arthur's
multiplicity formula for the packet $\mathcal{L} (\rho)$ produces.
$m_{\rm global} (\pi)$ is defined as the sum of $m (\mathcal{L}
(\pi'))$ where $\mathcal{L} (\pi')$ runs through all different
$L$-packets $\mathcal{L} (\pi')$ where $\pi$ and $\pi'$ are nearly
equivalent, i.e., $\pi'_{v} \cong \pi_{v}$ for almost $v$. This is
the multiplicity that occurs in the Arthur's multiplicity formula.
We will be mainly concerned with parameters which are {\it
tempered}, and for an explication of the Arthur's multiplicity
formula in this case we refer the readers to Lapid's paper
\cite{Lapid99}. We say that $G$ satisfies the multiplicity one if $m
(\pi) = m_{\rm global} (\pi) = 1$ for all cuspidal $\pi$. We know
the multiplicity one for $GL (n)$ (\cite{JPSS83}, \cite{JS81},
\cite{JS90}) and $SL (2)$ (\cite{Ra2003}).

\medskip

\subsection*{Some Heuristics}

Let us consider some heuristics. Let $F$ be a global field with the
conjectured Langlands group $\mathcal{L}_{F}$ (\cite{Arthur05},
\cite{Langlands2002}, \cite{Cl}, \cite{Ra}), and for each place $v$,
let $\mathcal{L}_{F_{v}}$ be the corresponding local group, which
can be taken to be the Weil group $\mathcal{W}_{F_{v}}$ for $v$
archimedean, and $\mathcal{W}_{F_{v}} \times SL_{2} (\mathbb{C})$
for $v$ non-archimedean. There should be natural morphism $j_{v}$:
$\mathcal{L}_{F_{v}} \to \mathcal{L}_{F}$ such that any global
parameter, i.e., a homomorphism $\phi$: $\mathcal{L}_{F} \to
{}^{L}G$ will induce local parameters $\phi_{v} = \phi \circ j_{v}$.

\medskip

Let $i$: ${}^{L}H \to {}^{L}G$ be an $L$-homomorphism such that
${}^{L}H^{\circ} (\mathbb{C})$ and ${}^{L}G^{\circ} (\mathbb{C})$
are connected, and the image of ${}^{L}H^{\circ} = \hat{H}$ is not
contained in any proper parabolic subgroup of ${}^{L}G$. Let $\pi$
be a cuspidal automorphic representation of $G (\mathbb{A}_{F})$
whose parameter $\phi_{\pi}$: $\mathcal{L}_{F} \to {}^{L}G$
satisfies the following:

\medskip

(a) $\phi_{\pi}$ factors through ${}^{L}H$.

\medskip

(b) ${\rm Im} (\phi_{\pi})$ is dense in ${}^{L}H$.

\medskip

This implies in particular that, for each place $v$, the local
parameter $\phi_{\pi_{v}}$: $\mathcal{L}_{F_{v}} \to {}^{L}G$
factors through ${}^{L}H$. Note that the Langlands group $\tilde{H}
(\pi)$ is essentially the Zariski closure of the image of ${}^{L}
H$.

\medskip

When $G = SO (2 N)$, ${}^{L}G^{\circ} = \hat{G} = SO (2 N,
\mathbb{C})$, and we can produce examples $\pi$ with multiplicities
$> 1$, if we assume functoriality for $i$: ${}^{L}H \to {}^{L}G$.
The main ideal is that when we have two $L$-parameters $\phi$ and
$\phi'$ which are locally conjugate but not globally conjugate in
image, they will produce such an example.

\medskip

To avoid the problem of the existence of $\mathcal{L}_{F}$ with
desired properties, we will restrict our attention to those $\pi$
whose parameters are tempered, and moreover naturally
representations of the absolute Galois (or Weil) group. This is
reasonable because our ultimate aim is to get nontrivial examples
which can be analyzed concretely.

\bigskip

\section{Potential Failure of Multiplicity One for Cusp Forms on $SO (2 N)$}
\label{S:6}

Let $G$ be a reductive algebraic group and $F$ an number field.
Again, for each field $k$, denotes the absolute Galois group as
$\mathcal{G}_{k}$. Let $H$ be a semisimple algebraic group such that
$\hat{H}$ embeds into $\hat{G}$. Let $i$ and $i'$ be two algebraic,
injective homomorphisms from $\hat{H}$ into $\hat{G}$ such that they
are locally conjugate but not globally conjugate in image. We also
denote by $i$ and $i'$ the induced $L$--homomorphisms from ${}^{L}H$
to ${}^{L}G$.

\medskip

We assume that:

\medskip

(1) The Langlands Functoriality for $i, i'$ holds, i.e., the
functorial transfer of automorphic forms for $H (\mathbb{A}_{F})$ to
$G (\mathbb{A}_{F})$ exists corresponding to the $L$-homomorphisms
$i, i'$: ${}^{L}H \to {}^{L}G$, at least for global parameters
$\phi$ attached to Galois representations.

\medskip

(2) The Arthur's conjecture, including the Arthur's multiplicity
formula for $G$ holds. i.e., the (global) multiplicity of each
cuspidal automorphic representation is exactly the same as given in
the Arthur's multiplicity formula. More precisely, to each
$L$-parameter $\rho$ of $G = SO (2 N)$, we associate a global
$L$-packets $\rho$ and a subrepresentation $L^{2} (\rho)$ of
$L^{2}_{\rm cusp}$ which is a direct sum of representations of
$\rho$ with multiplicities described by Arthur's multiplicity
formula. Moreover, if two parameters $\rho$ and $\rho'$ are not
conjugate in image, then $L^{2} (\rho)$ and $L^{2} (\rho')$ must be
orthogonal to each other.

\medskip

(3) There are two cuspidal automorphic representation $\pi$ and
$\pi'$ of $G (\mathbb{A}_{F})$, furnished by (1) from $\pi_{0}$ of
$H (\mathbb{A}_{F})$ with respective parameters $\phi$ and $\phi'$,
such that the image of $\phi$ is dense in $i (\hat{H})$.

\bigskip

Let $\phi_{0}$ be the parameter of $\pi_{0}$, and then we have $\phi
= i \circ \phi_{0}$ and $\phi' = i' \circ \phi_{0}$.

\medskip

First, $\pi$ and $\pi'$ are nearly equivalent. In fact, as $i$ and
$i'$ are locally conjugate, so are $\phi = \phi_{\pi} = i \circ
\phi_{0}$ and $\phi' = \phi_{\pi'} = i' \circ \phi_{0}$. Moreover,
$L^{S} (s, \pi, r) = L^{S} (s, \pi', r)$ for some finite set of
places of $F$.

\medskip

Moreover, $\pi$ and $\pi'$ cannot occur in the same constituent as
their parameters $\phi$ and $\phi'$ are not globally conjugate in
image. In fact, as the Zariski closures of $\phi$ and $\phi'$ are $i
(\hat{H})$ and $i' (\hat{H})$ respectively, they are not globally
conjugate in $\hat{G}$, and hence $L^{2} (\phi)$ and $L^{2} (\phi')$
are orthogonal to each other in $L^{2}_{\rm cusp} (G (F) \backslash
G (\mathbb{A}))$.

\medskip

So the global multiplicities of $\pi$ is at least $2$, and this
gives rise to an example to fail the multiplicity one for $G$.

\bigskip

Before we formulate a theorem in Galois version. we will start from
a global $l$--adic homomorphism $\phi_{0}: \mathcal{G}_{F} \to
\hat{H} (\bar{\mathbb{Q}}_{l}) \hookrightarrow {}^{L}H
(\bar{\mathbb{Q}}_{l})$ where $\hat{H}$ is the dual group of $H$,
and is always viewed as an algebraic subgroup of ${}^{L}H$ with
$\hat{H} (\mathbb{C}) = {}^{L}H^{\circ} (\mathbb{C})$.

\medskip

Put $\phi = i \circ \phi_{0}$ and $\phi' = i' \circ \phi_{0}$. Given
$\phi_{0}$ and any place $v$, we have a natural homomorphism from
$\mathcal{G}_{F_{v}}$ to $\mathcal{G}_{F}$, and by the restricting
$\phi$ defines $\phi_{v}: \mathcal{G}_{F_{v}} \to {}^{L}G
(\bar{\mathbb{Q}}_{l})$, and similarly $\phi'_{v}:
\mathcal{G}_{F_{v}} \to {}^{L}G (\bar{\mathbb{Q}}_{l})$.

\bigskip

The local parameter $\phi_{v}$: $\mathcal{G}_{F_{v}} \to {}^{L}G
(\bar{\mathbb{Q}}_{l})$ defines a homomorphism, again denoted
$\phi_{v}$ by abuse of notation, from $\mathcal{W}'_{F_{v}}$ to
${}^{L}G (\mathbb{C})$ (See Tate's article in Colvallis
\cite{Tate}). Hence $\phi_{v}$ should defines by functoriality to an
irreducible admissible representation $\pi_{v}$ of $G (F_{v})$. This
correspondence is known when $\phi_{v}$ is unramified. Same thing
works for $\phi'_{v}$.

\bigskip

\begin{theorem} \label{T:502}
Let $F$ be a global field and $l$ a prime not equal to the
characteristic of $F$. Let $H$ and $G$ be two algebraic reductive
groups, $i$ and $i'$ two algebraic injective morphism (or with
finite kernel) from $\hat{H} (\bar{\mathbb{Q}}_{l})$ to $\hat{G}
(\bar{\mathbb{Q}}_{l})$, and $\rho_{1}$ an $l$--adic Galois
homomorphism from $\mathcal{G}_{F}$ to $\hat{H}
(\bar{\mathbb{Q}}_{l})$, with all the images are semisimple. Assume
that,

\medskip

\textnormal{(1)} $i$ and $i'$ are locally conjugate but not globally
conjugate in image.

\medskip

\textnormal{(2)} The image of $\rho_{1}$ is Zariski dense in
$\hat{H}$.

\medskip

Let $\phi = i \circ \rho_{1}$ and $\phi' = i' \circ \rho_{1}$. Then
$\phi$ are $\phi'$ are locally conjugate but not globally conjugate
in image.

\medskip

Moreover, if $\phi$ and $\phi'$ are modular, i.e., if they are
associated to $\pi$ and $\pi'$, then $\pi$ and $\pi'$ are nearly
equivalent, and moreover give rise to multiplicity $\geq 2$ in the
space of cusp forms on $G (\mathbb{A}_{F})$.
\end{theorem}

\bigskip

The proof is similar to the discussion above.

\medskip

\emph{Proof}:

\medskip

First, $\phi$ and $\phi'$ are locally conjugate. In fact, as $i$ and
$i'$ are locally conjugate, and all the images $\rho_{1} (g)$ are
semisimple, then $\phi = i \circ \rho_{1}$ and $\phi' = i' \circ
\rho_{1}$ are also locally conjugate.

\medskip

Moreover, $\phi$ and $\phi'$ are not globally conjugate in image. In
fact, the Zariski closures of $\phi$ and $\phi'$ are $i (\hat{H}
(\bar{\mathbb{Q}}_{l}))$ and $i' (\hat{H} (\bar{\mathbb{Q}}_{l}))$
respectively, and they are not globally conjugate in image in
$\hat{G}$.

\medskip

Finally, if $\phi$ and $\phi'$ are modular, and are associated to
cuspidal automorphic representations $\pi$ and $\pi'$ of $G
(\mathbb{A}_{F})$ respectively, then $\pi$ and $\pi'$ are locally
conjugate and hence nearly equivalent. However, they give give rise
to multiplicity $\geq 2$ in the space of cusp forms on $G
(\mathbb{A}_{F})$ as the images of $\phi$ and $\phi'$ are not
conjugate.

\qed \bigskip

In particular, if $(\hat{H}, \hat{G})$ comes from our list from
Theorem ~\ref{TM:A}, i.e, $H$ is the simply connected algebraic Lie
group such that $\hat{H}$ is one of the simple adjoint Lie group
with Lie algebra $A_{4 n}\ (n \ge 1)$, $B_{2 n}\ (n \ge 2)$, $C_{2
n}\ (n \ge 2)$, $E_{6}$, $E_{8}$, $F_{4}$ and $G_{2}$, $i$ is the
adjoint representation and $i' = C \circ i$ where $C$ is some outer
automorphism of $G = SO (\mathfrak{g}, \kappa) = SO (2 N)$, which is
also the conjugation by some $g_{0} \in O (2 N) - SO (2 N)$.
Applying Theorem ~\ref{TM:A}, $i$ and $i'$ are locally conjugate but
not globally conjugate. Thus Theorem ~\ref{T:502} and hence Theorem
~\ref{TM:C} follow. Note that, if we assume the global Langlands for
$G = SO (2 N)$ here, or equivalently, first, assume the global
Langlands for $GL (2 N)$, then assume the Langlands functoriality
for the descent from $SO (2 N)$ to $GL (2 N)$, then Theorem
~\ref{TM:C} will also apply automatically.

\bigskip

\section{An example: Case $B_{2}$}

Now we come to the case when $\mathfrak{g} = B_{2} = C_{2}$ where
$\mathfrak{g} = so (5) = sp (4)$, $H = Sp (4)$ with $\hat{H} = SO
(5)$. As ${\rm dim}(\mathfrak{g}) = 10$, the adjoint representation
$i$ is also equivalent to the exterior square $\Lambda^{2}$. Then
Theorem ~\ref{T:502} applies. So the assumptions we need in this
case for the theorem are:

\medskip

(I) The functoriality from $Sp (4)$ to $SO (10)$ holds for the
exterior square or the adjoint representation $i$ from ${}^{L}Sp (4)
= SO (5) \to {}^{L}SO (10) = SO (10)$.

\medskip

(II) The functoriality from $Sp (4)$ to $SO (10)$ holds also for $i'
=\tau \circ i$ for some outer Lie automorphism $\tau$ of $SO (10)$.

\medskip

(II) There is an $l$--adic Galois homomorphism $\rho_{1}$:
$\mathcal{G}_{F} \to SO (5, \bar{\mathbb{Q}}_{l})$ whose image is
dense.

\bigskip

As $SO (5) \cong PSp (4) = Sp (4)/\{\,\pm1\,\}
= GSp (4)/Z$, so (II) can be replaced by:

(II$'$) there is a $4$--dimensional $l$--adic Galois representation
$\rho_{1}$ such that, the Zariski closure of its image contains $Sp
(4)$.

\bigskip

Let's quote the following lemma which is an easy result in the Lie
representation theory.

\begin{lemma}  \label{T:702}
Let $k$ be an algebraically closed field, and $\rho$ be a
$4$-dimensional irreducible continuous representation of a connected
semisimple algebraic group $H (k)$ over $k$ with finite kernel. Then
up to equivalence, the pair $(H, \rho)$ is one of the following:

\medskip

{\rm (1)} $H = SL (4)$, and $\rho$ is the standard representation.

\medskip

{\rm (2)} $H = SL (2)$, and $\rho = {\rm sym}^{3}$.

\medskip

{\rm (3)} $H = Sp (4)$, and $\rho$ is the standard representation.

\medskip

{\rm (4)} $H = SL (2) \times SL (2)$, and $\rho = \rho_{0} \otimes
\rho_{0}$ where $\rho_{0}$ is the standard representation of $SL
(2)$.

\medskip

{\rm (5)} $H = SO (4)$, and $\rho$ is the standard representation.

\medskip

Moreover, when $\rho$ is symplectic, only {\rm (2)} and {\rm (3)}
are possible, and when $\rho$ is orthogonal, only {\rm (4)} and
{\rm (5)} are possible.
\end{lemma}

\qed \bigskip

\begin{proposition} \label{T:505}

\medskip

Let $\rho_{1}$ be an irreducible continuous $4$-dimensional
representation of $\mathcal{G}_{F}$ over $\bar{\mathbb{Q}}_{l}$, and
assume the following:

\medskip

\textnormal{(a)} The Zariski closure of ${\rm Im} (\rho_{1})$ is a
reductive algebraic group.

\medskip

\textnormal{(b)} $\rho_{1}$ is (essentially) self dual, and is of
$GSp (4)$ type;

\medskip

\textnormal{(c)} $\rho_{1}$ is not twist equivalent to ${\rm
sym}^{3} (\rho_{2})$ for any continuous representation $\rho_{2}$ of
$\mathcal{G}_{F}$ over $\bar{\mathbb{Q}}_{l}$.

\medskip

Then the Zariski closure of ${\rm Im} (\rho_{1})$ contains $Sp (4)$.

\end{proposition}

\bigskip

\emph{Proof}: This can be deduced from the results of \cite{Ra2003}.
Let us give the argument for completeness.

\medskip

Let $H$ be the Zariski closure of ${\rm Im} (\rho_{1})$, and $\rho$
be the embedding of $H$ into $GL (4)$. Moreover, denote $H' =
H^{ss}$ be a semisimple part of $H$. Note that $H'
(\bar{\mathbb{Q}}_{l})$ is unique up to conjugacy.

\medskip

Then $\rho$ restricted to $H' = H^{ss}$ is also irreducible. Hence,
from Lemma ~\ref{T:702}, $H' (\bar{\mathbb{Q}}_{l})$ has to be
isomorphic to $SL (4)$, $Sp (4)$, $SO (4)$ or $SL (2) \times SL
(2)$. Note that (b) rules out $SL (4)$, $SO (4)$ and $SL (2) \times
SL (2)$, and (c) rules out $SL (2)$. Thus only $Sp (4)$ remains, and
hence the proposition.

\qed \bigskip

\emph{\textbf{Remark}}: Our goal in a sequel will be to first
construct modular Galois representations satisfying the Assumptions
(1), (2) of Theorem ~\ref{T:502} in the case $B_2$, which we think
is possible over $\mathbb{Q}$ and over function fields, by  starting
with cusp forms with suitable discrete series components, and then
specialize to the function field case where we can appeal to the
deep work of Lafforgue (\cite{Laf}) on $GL(n)$ and also use the
recent striking results on the automorphic transfer to $GL(n)$ from
classical groups of Cogdell, Kim, Piatetski-Shapiro and Shahidi
(\cite{CoKPSS}), and also the backward lifting (automorphic descent)
of Ginzburg, Rallis and Soudry (for the survey, see \cite{S2005}).

\end{document}